\newif\ifsmfart
\IfFileExists{smfart.cls}
   {\documentclass[12pt,english]{smfart}
    \IfFileExists{smfenum.sty}{\usepackage{smfenum}}{}
    \usepackage{bull}
    \smfarttrue}
{\message{^^J*** It would be better to typeset this file with smfart.cls ***^^J^^J}

\documentclass[12pt]{amsart}}

\usepackage{amssymb}
\usepackage{epsfig}
\usepackage{graphicx}

\numberwithin{equation}{section}

\input xy
\xyoption{all}

\advance\headheight 2pt
\calclayout
\allowdisplaybreaks[3]

\theoremstyle{plain}
\newtheorem{prop}[subsection]{Proposition}

\newtheorem{thm}[subsection]{Theorem}
\newtheorem{coro}[subsection]{Corollary}

\newtheorem{lemm}[subsection]{Lemma}

\newtheorem{defn}[subsection]{Definition}

\theoremstyle{definition}
\newtheorem{prob}[subsection]{Problem}
\newtheorem{ques}[subsection]{Question}
\newtheorem{conj}[subsection]{Conjecture}

\theoremstyle{remark}
\newtheorem{rem}[subsection]{Remark}

\newtheorem{exam}[subsection]{Example}

\def\End{{\rm End}}
\def\rk{{\rm rk\,}}
\def\Tr{{\rm Tr}}

\newcommand{\la}{\lambda}

\newcommand{\C}{\Bbb C}
\newcommand{\Q}{\Bbb Q}
\newcommand{\Z}{\Bbb Z}

\def\Fr{{\rm Fr}}

\def\Sym{{\rm Sym}}
\def\Gal{{\rm Gal}}

\def\char{{\rm char}}

\def\ra{\rightarrow}

\def\C{{\mathbb C}}

\def\P{{\mathbb P}}
\def\Q{{\mathbb Q}}

\def\Z{{\mathbb Z}}
\def\C{{\mathbb C}}
\def\N{{\mathbb N}}

\def\Aut{{\rm Aut}}
\def\Pic{{\rm Pic}}

\def\ord{{\rm ord}}

\makeatother
\makeatletter

\def\Title    {Rational curves and points on K3 surfaces}
\def\Author   {Fedor Bogomolov and Yuri Tschinkel}
\def\Subject  {Algebraic geometry, number theory}
\def\Keywords {Curves, algebraic points, K3 surfaces}

\newif\ifpdf
   \ifx\pdfoutput\undefined
   \pdffalse
   \else
           \pdfoutput=1
            \pdftrue
   \fi
      \ifpdf
\usepackage[pdftex,
hyperindex=true,
pdftitle={\Title},pdfauthor={\Author},pdfsubject={\Subject},
pdfkeywords={\Keywords},pdfpagemode={UseOutlines},
bookmarks=true,bookmarksopen=true,
bookmarksopenlevel=0,bookmarksnumbered=true,
pdfstartview={FitV},
colorlinks=true
]{hyperref}
 \else
 \usepackage{hyperref}
 \fi

\usepackage{graphicx}

\author{Fedor Bogomolov}
\address{Courant Institute of Mathematical Sciences, N.Y.U. \\
 251 Mercer str. \\
 New York, NY 10012, U.S.A.}
\email{bogomolo@cims.nyu.edu}

\author{Yuri Tschinkel}
\address{Mathematisches Institut \\
         Universit\"at G\"ottingen\\
         Bunsenstr. 3-5\\
         37073 G\"ottingen, Germany}
\email{yuri@uni-math.gwdg.de}

\title[Rational curves and points on K3 surfaces]
{Rational curves and points \\ on K3 surfaces}

\begin{document}

\date{\today}



\begin{abstract}
We study the distribution of algebraic points on K3 surfaces.
\end{abstract}

\maketitle

\tableofcontents

\setcounter{section}{0}
\section{Introduction}
\label{sect:introduction}

Let $k$ be a field and $\bar{k}$  
a fixed algebraic closure of $k$.
We are interested in connections 
between geometric properties
of algebraic varieties and their
arithmetic properties over $k$, 
over its finite extensions $k'/k$ or over $\bar{k}$.
Here we study certain varieties 
of intermediate type, namely K3 surfaces and their 
higher dimensional generalizations, Calabi-Yau varieties.   

\

To motivate the following discussion, consider a K3 surface $S$
defined over $k$.  In positive characteristic, 
$S$ may be unirational and covered by rational 
curves. 
If $k$ has characteristic zero, 
then $S$ contains countably many rational 
curves, at most finitely many in each homology class of $S$ 
(the counting of which is an interesting problem in 
enumerative geometry, see \cite{beauville}, \cite{bryan}, \cite{chen}, 
\cite{yau}). Over uncountable fields, there may, of course, exist $k$-rational 
points on $S$ not contained in any rational curve defined over $\bar{k}$. 
The following extremal statement, proposed by the first author 
in 1981, is however still a logical possibility:

\

{\it Let $k$ be either a 
finite field or a number field.
Let $S$ be a K3 surface defined over $k$.
Then every $\bar{k}$-rational point on $S$ lies on  
some rational curve $C\subset S$, defined over $\bar{k}$.}

\

In this note we collect several representative examples 
illustrating this statement. One of our results is:

\begin{thm}
\label{thm:main}
Let $S$ be a Kummer surface over a finite field $k$.
Then every $s\in S(\bar{k})$ lies on a rational curve $C\subset S$, 
defined over $\bar{k}$. 
\end{thm}

Using this theorem we produce 
examples of non-uniruled but ``rationally connected'' surfaces over finite fields
(any two algebraic points can be joined by a chain of rational curves).

\

\noindent
{\bf Acknowledgments:} 
We are grateful to Brendan Hassett, Ching-Li Chai, Nick Katz, Barry Mazur and 
Bjorn Poonen for their interest and comments.

\section{Preliminaries: abelian varieties}
\label{sect:abel}

In this section we collect some facts concerning abelian 
varieties. Our basic reference is \cite{mumford}. 

\

Let $A$ be an abelian variety over $\bar{k}$, 
$A[n]\subset A(\bar{k})$
the set of the $n$-torsion points of $A$.
If $k$ is finite, then every point in $A(\bar{k})$ 
is a torsion point.
For every torsion point $x\in A(\bar{k})$ let 
$$
\ord(x):=\min\{ n\in \N\,|\, nx=1\}
$$
be the order of $x$. 
An elliptic curve $E$ over a field $k$ of characteristic $p$ is called
{\em supersingular} if its $p$-rank is zero, an  
abelian variety over $k$ is called supersingular
if it is $\bar{k}$-isogenous to a product of supersingular elliptic curves.

Recall that every abelian variety $A$, defined over $\bar{k}$, 
is isogenous to a product of {\em simple} abelian varieties (over $\bar{k}$). 
The ring $\End_{\bar{k}}(A)$ of $\bar{k}$-endomorphisms 
of a simple abelian variety $A$ is a maximal order in some 
finite dimensional division algebra $D$ over $\Q$, 
of dimension $d^2$ over its 
center $Z_D$. Here $Z_D$ is a finite extension of $\Q$ of degree $e$.   
If $k$ is finite, then $ed=2\dim(A)$, that is, 
$A$ has complex multiplication.

\begin{rem}
\label{rem:oort}
In our applications, we will use hyperelliptic curves contained
in abelian varieties. Recall that, over an algebraically closed field,  
every abelian surface is the Jacobian of a (possibly reducible) hyperelliptic curve 
(see \cite{weil}). 
A simple argument shows that for any abelian variety $A$ 
over $\mathbb C$ (in any dimension) 
the set of hyperelliptic 
curves which have a Weierstrass point (fixed point of the hyperelliptic involution)
at the origin of $A$ is at most discrete. Indeed, a family of such curves
would give a uniruled surface in the quotient $A/\sigma$
(under the standard multiplication by $-1$ map $\sigma$). 
However, $A/\sigma$ has nontrivial holomorphic two-forms 
which extend to its desingularization. 
They cannot restrict to a uniruled surface. 
A much more precise result is proved in \cite{pirola}:
a generic principally polarized abelian variety of dimension $\ge 3$  
over $\C$
does not contain {\em any} hyperelliptic curves. Similarly, 
a generic abelian variety of dimension $\ge 4$ does not contain
trigonal curves \cite{alz-pirola}.
A similar result holds in positive characteristic, over {\em large} fields, like 
an algebraic  closure of $\bar{\mathbb F}_q(t)$
\cite{oort-dejong}. 

It could still be possible that over an algebraic closure of a 
{\em finite} field, every abelian variety contains a hyperelliptic curve.
\end{rem}

\

Let $C$ be a smooth projective geometrically connected curve
of genus $\mathsf g=\mathsf g(C)$ over a field $k$ and $J=J_C$ 
the Jacobian of $C$. Throughout, we assume that $C(k)\neq \emptyset$ and choose
a point $c_0\in C(k)$ which we use to identify the degree $n$ Jacobian $J_C^{(n)}$ with $J_C$
and to embed $C\ra J_C$. 
Consider the maps 

\

\centerline{
\xymatrix{
C^n \ar[r]              & C\times \Sym^{(n-1)}(C) \ar[r]^{\,\,\,\,\,\,\,\phi_n} & 
\Sym^{(n)}(C) \ar[r]^{\,\,\,\,\,\,\,\varphi_n} & J_C^{(n)} \\
c=(c_1,\ldots, c_n) \ar[r] & (c_1,c_2+\cdots +c_n) \ar[r]   & (c_1+\cdots +c_n)\ar[r] & [c] 
}}

\

The map $\phi_n$ is a (finite, non-Galois for $n\ge 3$) cover of degree $n$.  
For all $n\ge 2\mathsf g+1$, the map $\varphi_n$ is a $\P^{n-\mathsf g}$-bundle and
the map $C^n\ra J^{(n)}(C)$ has geometrically irreducible fibers
(see \cite{katz1}, Cor. 9.1.4, for example).

\begin{lemm}
\label{lemm:main}
Let $k$ be a number field or a  (sufficiently large) finite field
with a fixed algebraic closure $\bar{k}$. 
Let $C$ be a smooth projective geometrically connected curve over $k$ 
of genus $\mathsf g=\mathsf g(C)\ge 1$,  with Jacobian $J=J_C$. 
For every  $n\ge 2\mathsf g+1$ and every point $x\in J(k)$ 
there exists a point $y\in \varphi_n^{-1}(x)$, defined over $k$, 
whose preimage $c=(c_1,\ldots, c_n)\in C^n(\bar{k})$ gives rise to 
a $k$-irreducible degree $n$ zero-cycle $c_1+\cdots +c_n$ on $C$. 
\end{lemm}

\begin{proof}
Let $x\in J(k)$ be a point and $\P_x=\varphi_n^{-1}(x)\subset \Sym^{(n)}(C)$ 
the fiber over $x$. The restriction $\phi_{n,x}$ of 
$\phi_n$ to $\P_x$ is a (nontrivial) 
cover of degree $n$. 
Using Chebotarev density 
(or an equidistribution theorem, as in \cite{katz1}, Th. 9.4.4)
we find that, if $k$ is sufficiently large (either a number field, or
a finite field with $\ge \mathsf c(\mathsf g)$ elements, 
where $\mathsf c(\mathsf g)$ is an 
explicit constant depending only on the genus $\mathsf g$), 
there exists a point $y\in \P_x$ such that the fiber $\phi_{n,x}^{-1}(y)$
is irreducible over $k$. 
\end{proof}

\begin{rem}
The same statement holds for quasi-projective curves $C$ and their
(generalized) Jacobians. 

A similar result (over finite fields) 
has recently been used in \cite{poonen}, Lemma 5.
\end{rem}

We reformulate the above theorem in a more convenient form:

\begin{coro}
\label{coro:K}
Let $C$ be a curve of genus $\mathsf g=\mathsf g(C)\ge 2$ 
over a number field $K$ and  
$J=J_C$ its Jacobian. 
Assume that $C$ has a zero-cycle of degree $n=2\mathsf g+1$  
over $K$ and use this cycle to identify $J^{(n)}=J$ and the embedding $C\ra J$.  
For any point $x\in J(K)$ there exists an extension $K'/K$
of degree $n$ and a point $c\in C(K')$ such that 
the cycle $\Tr_{K'/K}(c)=x\in J^{(n)}(K)=J(K)$. 
\end{coro}

\begin{coro}
\label{coro:need}
Let $C$ be a curve of genus $\mathsf g(C)\ge 2$ over a
(sufficiently large) finite field $k$,  
$J=J_C$ its Jacobian and $x\in J(k)$ a point. 
Choose a point $c_0$ on $C(k)$ and use it to identify 
$J_C=J^{(n)}_C$, for all $n$, and to embed $C\ra J_C$.
Then there exist 
a point $c\in C(\bar{k})$ and infinitely many endomorphisms 
$\Phi\in \End_{\bar{k}}(J)$ such that $\Phi(c)=x$.  
\end{coro}

\begin{proof}
For any $n \ge  2{\mathsf g}(C)+1$ consider the surjective map $\varphi_n$.
Let $x \in J^{(n)}_C({k})$  be a point and
$\P_x=\varphi_n^{-1}(x)$ the projective space over $x$, 
parametrizing all zero-cycles equivalent to $x$.
Extending $k$, if necessary, we find a $y\in \P_x(k)$
such that the fiber 
$$
\phi_n^{-1}(y)= \{ (c_1, c_2+\cdots + c_n), (c_2, c_1+\cdots + c_n), 
\ldots, (c_n, c_1+\cdots + c_{n-1})\} 
$$
is irreducible over $k$ (see \cite{katz1}, Th. 9.4.4).

We can write $y=\sum_{g\in G} c^{g}$, with  
$c:=c_1\in C(k')$, where $k'/k$ is the (unique)
extension of $k$ of degree $n$ and
$g$ runs through the (cyclic)
Galois group $G:=\Gal(k'/k)$ over $k$.  Note that $g$ are 
simply (distinct) powers of the Frobenius morphism. 
We have 
$$
y=\sum_{j=0}^{n-1} {\rm Fr}^j(c).
$$ 
Now note that the Frobenius morphism ``lifts'' to a $k'$-endomorphism of $J$, that
is, the Frobenius endomorphism $\tilde{\rm Fr}\in \End_{\bar{k}}(J)$ acts on $J(k')$
in the same way as the Galois automorphism ${\rm Fr}\in \Gal(k'/k)$. 
Put 
$$
\Phi:=\sum_{j=0}^{n-1} \tilde{\rm Fr}^j,
$$ 
as an element of $\End_{\bar{k}}(J)$.
Changing $n$, we get infinitely many such endomorphisms (over $\bar{k}$). 
\end{proof}

\begin{rem}
In particular, Corollary~\ref{coro:need} implies that 
if $\ord(x)=m$ then,  for any embedding of $C \subset J_C$, there
exist infinitely many points in $C(\bar{k}) \subset J_C(\bar{k})$ 
whose order is divisible by $m$. Indeed, 
notice that $\ord(c) = \ord(c^{g})$, for all $g \in \Gal(k'/k)$, 
(with respect to some group law on the set $J(k)$). 
Since the order of $x$ (and $y$) is $m$ the order  
$\ord(c)$ is divisible by $m$.

A related result has been proved
in \cite{anderson}:
Let $\ell$ be a prime, $C$ a curve, $J=J_C$ its Jacobian  
(defined over a finite field $k$), $C\subset J$ a {\em fixed} embedding
and $\la\,:\, J(\bar{k})\ra J(\bar{k})_{\ell}$ the projection 
onto the $\ell$-primary part. Then the map 
$\la\,:\, C(\bar{k})\ra J(\bar{k})_{\ell}$ is surjective.  
It was noticed in \cite{pop} that this fact implies that
any positive-dimensional subvariety of a geometrically simple abelian variety
(over a finite field) contains infinitely many points 
of pairwise prime orders. 
\end{rem}

The same argument gives a statement very much in the spirit of \cite{katz}:

\begin{coro}
\label{coro:katz}
Let $C$ be a curve of genus $\mathsf g$ over a sufficiently 
large finite field $k$, $J=J_C$ its Jacobian and  
$k'/k$ the (unique) degree $2\mathsf g+1$ extension of $k$.
Then there exists a morphism $C\rightarrow J$ (depending on $k$)
such that $J(k)\subset C(k')$.
\end{coro}

\section{Preliminaries: K3 surfaces}

We assume that the characteristic of $k$ is either zero or at least  $7$.

\begin{defn}
A connected simply-connected projective algebraic surface with trivial 
canonical class is called a K3 surface. 
A K3 surface $S$ with $\rk \Pic(S)=22$ 
is called {\em supersingular}. 
\end{defn}

\begin{exam}
\label{exam:typical}
Typical K3 surfaces are 
double covers of $\P^2$ ramified in a smooth curve of degree 6, 
smooth quartic hypersurfaces in $\P^3$ or
smooth intersections of 3 quadrics in $\P^5$.

Another interesting series of examples is given by 
(generalized) Kummer surfaces: desingularizations of 
quotients of abelian surfaces by certain finite group actions 
(see Proposition~\ref{prop:katsu}). 
\end{exam}

\begin{rem}
If $S$ is a K3 surface over a field of characteristic zero, 
then $\rk\Pic(S)\le 20$. 
An example of a  
supersingular $S$ over a field of positive characteristic
is given by a desingularization of $A/\sigma$, where $A$ is
a supersingular abelian variety and $\sigma$ the standard involution
(multiplication by $-1$ map). 
\end{rem}

\begin{rem}
\label{rem:pos}
If $\rk\Pic(S)<22$ (and hence $\le 20$) then 
the Brauer group of $S$ has 
nontrivial transcendental part. In particular, 
$S$ is not uniruled. This is always the case in 
characteristic zero.  
Over fields of positive characteristic, there may  
exist uniruled K3 surfaces, with necessarily 
the maximal possible Picard number $\rk\Pic(S)=22$ 
(all cycles are algebraic); and they are therefore 
supersingular \cite{shaf-rudakov1}, \cite{artin}.

In characteristic 2, every supersingular K3 surface is 
unirational \cite{shaf-rudakov}. It is conjectured that
all supersingular K3 surfaces are unirational.
A generalized Kummer surface $S \sim A/G$
is uniruled iff it is unirational  
iff the corresponding abelian surface $A$ is 
supersingular \cite{shioda}, \cite{katsura}. 
\end{rem}

\

\section{Construction}
\label{sect:const}

We recall the classical construction of 
special K3 surfaces, called {\em Kummer} surfaces. Let 
$A$ be an abelian surface,
$$
\begin{array}{rcc}
\sigma\,:\, A & \ra &  A\\
           a & \mapsto & -a 
\end{array}
$$
the standard involution. 
The set of fixed points of $\sigma$ 
is exactly $A[2]$. 
The blowup $S:=\widehat{A/\sigma}$ of 
the image of $A[2]$  in the quotient $A/\sigma$ 
is a smooth K3 surface $S$, called a {\em Kummer} surface:
$$  
\pi\,:\, A/\sigma  \ra  S,  \hskip 1cm 
\hat{\pi}\,:\, \widehat{A /\sigma}   \ra S.
$$ 

\

Consider rational curves in $A /\sigma$.

\begin{lemm} 
\label{lemm:ratt}
Rational curves $C\in A/\sigma$ correspond to
hyperelliptic curves $\tilde{C} \subset A$ 
containing a two-torsion point $P\in A[2]$.
\end{lemm}

\begin{proof} 
The hyperelliptic involution on $C$ acts as an involution 
$\sigma\,:\, x \to -x$ on the Jacobian $J_C$ and hence also on the 
abelian subvariety which is the image
of $J_C$ in $A$. In particular, the involution $\sigma$ on $A$
induces the standard hyperelliptic involution on $C$. Hence 
$C /\sigma = \P^1$ is rational and defines a rational curve in $ A/\sigma$.
Conversely, if $\P^1\in A/\sigma$ is rational then the preimage of  $\P^1$
in $A$ is irreducible (since $A$ doesn't contain
rational curves). Thus $\P^1 = \tilde{C}/\sigma$ and $\tilde{C}$ 
is hyperelliptic and all ramification points of the map 
$\tilde{C} \to C=\P^1$ are contained among the two-torsion points 
$A[2]\cap \tilde{C}$.
\end{proof}

\begin{thm}
\label{thm:k}
Let $k$ be a finite field, $C$ a curve of genus $2$ defined
over $k$, $J=J_C$ its Jacobian surface and $S\sim J/\sigma$ 
the associated Kummer surface.
Then every algebraic point $s\in S(\bar{k})$
lies on some rational curve, defined over $\bar{k}$.
\end{thm}

\begin{proof}
Let $s\in S(\bar{k})$ be an algebraic point (on the complement to the
16  exceptional curves) and $x\in J(\bar{k})$ one of its preimages. 
We have proved in Corollary~\ref{coro:need} 
that for every $x\in  J(\bar{k})$ there
is an endomorphism $\Phi\in \End_{\bar{k}}(J)$ such that 
$\Phi \cdot C(\bar{k})$ contains $x$ (note that $\Phi$ commutes with
the involution $\sigma$). 
The image of the curve $\Phi\cdot C$ in $S$ contains $s$. 
\end{proof}

\begin{coro} 
Let $S$ be a Kummer surface over
a finite field $k$. 
There are infinitely many 
rational curves (defined over $\bar{k}$) 
through every point in $S(\bar{k})$ 
(in the complement to the 16 exceptional curves).  
If $S$ is non-uniruled, these curves 
don't form an algebraic family. 
\end{coro}

In addition to quotients $A/\sigma$, there exist (generalized) Kummer 
K3 surfaces obtained as desingularizations of abelian surfaces under
actions of other finite groups. Such actions
(including positive characteristic) 
have been classified:

\begin{prop}[see \cite{katsura}]
\label{prop:katsu}
Let $A$ be an abelian surface over a field $k$ and  
$G$ a finite group  acting on $A$ such that
the quotient $A/G$ is birational to a K3 surface.
If $\char(k)>0$ then $G$ is one of the following:
\begin{itemize}
\item a cyclic group of order $2,3,4,5,6,8,10,12$;
\item a binary dihedral group $(2,2,n)$ with $n=2,3,4,5,6$;
\item a binary tetrahedral group $(2,3,3)$;
\item a binary octahedral group $(2,3,4)$;
\item a binary icosahedral group $(2,3,5)$.
\end{itemize}
If $\char(k)=0$ then $G$ is one of the following:
\begin{itemize}
\item a cyclic group of order $2,3,4,6$;
\item a binary dihedral group $(2,2,n)$ with $n=2,3$;
\item a binary tetrahedral group $(2,3,3)$.
\end{itemize}
\end{prop}

The groups listed above do indeed occur. 
 
\begin{coro}
If $S\sim A/G$ is a generalized Kummer K3 surface  
over a finite field $k$ (of characteristic $\ge 7$) 
then every algebraic point on $S$ lies on 
infinitely many rational curves, 
defined over $\bar{k}$. 
\end{coro}

\begin{proof}
By Remark~\ref{rem:pos}, a supersingular Kummer K3 surface is uniruled
and the claim follows.  
By Lemma 6.2 in \cite{katsura} if $S$ is not supersingular and  
$G$ is divisible by $2$ then $G$ has a unique element of order two, acting 
as the standard involution. An argument as in the proof of 
Theorem~\ref{thm:k} applies to show
that every algebraic point lies on a rational curve.
The Kummer K3 with  $G=\Z/5$ are supersingular. 

It remains to consider $G=\Z/3$. 
In this case $A=E_0\times E_0$ with $E_0$ the elliptic curve $y^2=x^3-1$, 
with complex multiplication by $\Z[\sqrt{-3}]$. Let $C$ be the genus two 
curve given by $y^2=x^6-1$. Its Jacobian is (isogenous to) $E_0\times E_0$. 
The natural $\Z/3$-action has eigenvalues $\zeta_3$, $\zeta_3^2$ and the quotient of 
$C$ by this action is a rational curve.  
Applying the argument of Corollary~\ref{coro:need} and
endomorphisms (sums of powers of Frobenii, they commute with 
the $\Z/3$-action) we obtain our claim. 
\end{proof}

\begin{rem}
\label{rem:others}
There exist K3 surfaces of non-Kummer type, which are dominated by 
Kummer K3 surfaces. Clearly, these have the same property. 
\end{rem}

\begin{rem}
We don't know whether or not {\em every} algebraic K3 surface 
contains infinitely many rational curves 
(elliptic K3 surfaces do, see \cite{bt}).
Thus it is tempting to consider the question of 
lifting of rational curves on K3 surfaces
from characteristic $p$ to characteristic zero. This 
is analogous to deformations of rational curves 
on K3 surfaces over $\C$, where the answer is, roughly speaking, 
that the rational curve deforms as long as its homology class 
remains algebraic (see, for example,  
\cite{mori-mukai}, \cite{bt}). 
If this principle applies, then every K3 surface $S$  
over $\bar{\Q}$ which reduces to a 
Kummer K3 surface modulo at least one prime, has 
infinitely many rational curves. 

It is known that primitive classes in $\Pic(S)$
of a {\em general} K3 surface $S$ over $\C$
are represented by rational curves with at worst nodal singularities
(see \cite{yau}, \cite{chen}, for example). In particular, a general 
polarized $S$ with $\rk\Pic(S)\ge 2$ has infinitely many rational curves. 
See, however, \cite{ellenberg} for examples of surfaces
with $\rk\Pic(S_{\bar{\Q}})=1$.  
\end{rem}

\begin{rem}
\label{rem:fail}
Theorem~\ref{thm:k} fails if $k={\mathbb F}_q(t)$ and if $S$ is an isotrivial 
Kummer surface over $k$. 
Indeed, we can think of $S$ as a fibration over
$\P^1$ (over $\mathbb{F}_q$) and choose a rational curve $C_0$ 
in a (smooth) fiber $S_0$ of this fibration.  Then,
for some $C_0'$ over $C_0$, there is a surjective map
$C_0'\times \P^1\ra S_0$. 
However, one can choose a non-uniruled (non-supersingular) $S_0$.
\end{rem}

\section{Surfaces of general type}
\label{sect:gen-type}

Using similar ideas we can construct non-uniruled surfaces $S$
of general type over finite fields $k$ such that every algebraic point
$s\in S(\bar{k})$ lies on a rational curve and any two points can 
be connected by a chain of rational curves. (However, the degrees of these
curves cannot be bounded, {\em a priori}).

For simplicity, let us assume that $p:=\char(k)\ge 5$. 
Let $S_0$ be a unirational surface of general type over $k$, 
for example 
$$
x^{p+1}+y^{p+1}+z^{p+1}+t^{p+1}=0,
$$  
and $\P^2\ra S_0$ the corresponding 
(purely inseparable) covering of degree a power of $p$. 
Let $S_1$ be a non-supersingular (and therefore, non-uniruled) 
Kummer K3  surface admitting an abelian cover
onto $\P^2$ of degree prime to $p$ with Galois group $G$ 
(for example, a double cover). 

\begin{lemm}
\label{lemm:easy}
For any $n$ coprime to $p$, and any purely inseparable extension
$K/L$ we have a natural isomorphism 
$$
K^*/(K^*)^n =L^*/(L^*)^n.
$$
\end{lemm}
 
\begin{proof}
Indeed, there exists an $m\in \N$ such $L^*$ contains the 
$p^m$-powers of all elements of $K^*$. 
Since $p^m$ and $n$ are coprime the claimed isomorphism follows.  
\end{proof}

By Kummer theory, the extension of function fields 
$\bar{k}(S_1)$ over $K$ is obtained by 
adjoining roots of elements in $K^*$. 
The extension is defined modulo $(K^*)^n$, for some $n$ coprime to $p$. 
By Lemma~\ref{lemm:easy}, we can select a $\phi$ in $L^*:=\bar{k}(S_0)^*$, 
which gives this extension. 
In particular, we get a diagram
$$
\begin{array}{ccc}
S_1  &  \ra  &  S \\
\downarrow   &  & \downarrow \\
\P^2  & \ra &   S_0,
\end{array}
$$
where $S$ is a surface of general type 
(since the corresponding function field
is a separable abelian extension of degree coprime to $p$).
At the same time there is a surjective purely inseparable map $S_1\ra S$.
Surjectivity implies that there is a rational curve
(defined over $\bar{k}$) passing through every 
algebraic point of $S$ (to get {\em every} point we may need to pass to a blowup 
$\tilde{S}_1$ of $S_1$ resolving the indeterminacy of the dominant map $S_1\ra S$). 
Pure inseparability implies that 
$S$ is non-uniruled, since $S_1$ is non-uniruled.

\section{The case of number fields}
\label{sect:numer}

In this section, $K$ is a number field (with a fixed embedding into 
$\bar{\Q}$). For any nonarchimedean place  
$v$ of $K$ let $k_v$ be the residue field at $v$. 
Let $A$ be an absolutely simple abelian variety over $K$ of dimension $\mathsf g$. 
Assume that $\mathcal O:=\End_{\bar{\Q}}(A)$ is an 
order in a field $F$ with $[F : \Q] = 2\mathsf g$ (complex multiplication). 
Let $\mathsf S$ be a finite set of
places of good reduction of $A$,  
$A_v$ the corresponding abelian varieties over finite
fields $k_v$, for $v\in \mathsf S$.

\begin{lemm}
\label{lemm:appro}
Let $C$ be a curve of genus $\mathsf g(C)\ge 2$ defined over $K$ 
and $J=J_C$ its Jacobian. Assume that $J$ is absolutely simple and  
has complex multiplication.
Let $\{x_v\}_{v\in \mathsf S}$ be a set of smooth points $x_v\in J_v(k_v)$. 
Then there is an endomorphism $\Phi\in \End_{\bar{\Q}}(J)$ 
and a $c\in C(\bar{\Q})$ such that $\Phi(c)_v = x_v$ for any $v\in \mathsf S$.
\end{lemm}

\begin{coro}
In particular, let $C$ be a hyperelliptic curve over a number field $K$ with 
absolutely simple Jacobian $J=J_C$,
and let $\sigma$ the standard involution on $J$ (multiplication by $-1$), 
acting as a hyperelliptic involution on $C\subset J$. 
Assume that $J$ has complex multiplication.  
Then every point on $S\sim J/\sigma$ can be approximated
by points lying on rational curves, that is, for any finite set of
places $\mathsf S$ of good reduction and any set of points 
$s_v \in S(k_v)$ there is a point $s\in S(\bar{\Q})$ lying on 
a rational curve on $S$, defined over $\bar{\Q}$,  such that
$s$ reduces to $s_v$, for all $v\in \mathsf S$. 
\end{coro}

\begin{proof}[Proof of Lemma~\ref{lemm:appro}]
There is an $n\in \N$ (which we now fix) such that 
$x_v\in J_v[n]$, for all $v\in \mathsf S$.
Now we choose a sufficiently large number field $L/K$ 
so that $J(L)\ra J_{w}[n]$ is  a
surjection for all places  $w\in \mathsf S_L$ (the set of places 
of $L$ lying over $v\in \mathsf S$).
In particular, for all residue fields $k_w$, 
$J_w[n]\subset J_v(k_w)$ so that the corresponding Frobenii $\Fr_w$, considered as
(conjugacy classes of) elements of the Galois group $\Gal(\bar{L}/L)$,   
act trivially on $J_{w'}[n]$, for all $w'\in \mathsf S_L$. 

Consider the map 
$$
\psi_w\,:\, \End_{\bar{\Q}}(J) \to \End_{\bar{k}_w}(J_w).
$$ 
Since $J$ has complex multiplication, 
$\psi_w(\End_{\bar{\Q}}(J))$ contains 
the central subfield in 
$\End_{\bar{k}}(J_w)$.
Lift $\Fr_w\in \Gal(\bar{k}_w/k_w)$ to $\End_{\bar{k}_w}(J_w)$ and 
then (since it is a central element) to a global endomorphism
$\Psi_w\in \End_{\bar{\Q}}(J)$.

Now we put
$$
\Phi_w:=\sum_{i=0}^{m-1} \Psi_w^i \in \End_{\bar{\Q}}(J),
$$
with $m$ an integer $\ge 2\mathsf g+1$ and congruent to $1 \mod n$) 
and
$$
\Phi:=\prod_{w\in \mathsf S_L} \Phi_w
$$
(note that $\Phi_w$ commute in $\End_{\bar{\Q}}(J)$).

We have seen  in Corollary~\ref{coro:need} 
that for all $w$ and all $x_w\in J_w[n]$
there exists a point $c_w\in C(\bar{k}_w)$ such that
$\psi_w(\Phi_w)(c_w)=x_w$. Observe that 
in fact $\psi_w(\Phi)(c_w)=x_w$, for all $w\in \mathsf S_L$.
Indeed, since $\Phi=\prod_{w'\neq w} \Phi_{w'} \cdot \Phi_w$, 
and $\psi_w(\Phi_{w'})$ are trivial on 
$x_w=\Phi_{w}(c_w)\in J_w[n]$ 
(for all $w'\in \mathsf S_L$), we get the claim.

To conclude, choose a point $c\in C(\bar{\Q})$ reducing to $c_v\in C(\bar{k}_v)$.
Then $\Phi(c)$ reduces to $x_v$, for all $v\in \mathsf S$.
\end{proof}

\section{Higher dimensions}
\label{sect:higher}

\begin{thm}
\label{thm:kk}
Let $k$ be a finite field, $C$ a hyperelliptic curve of genus $\ge 2$ over $k$,
$J=J_C$ its Jacobian abelian variety, $\sigma$ the standard involution on $J$ and 
$S=J/\sigma$ the associated (generalized) Kummer surface. 
Then every rational point $s\in S(\bar{k})$
lies on some rational curve, defined over $\bar{k}$. 
\end{thm}

\begin{proof}
See Theorem~\ref{thm:k}.
\end{proof}

\begin{defn}
A smooth projective variety $V$ is called Calabi-Yau
if its canonical class is trivial and 
$h^0(\Omega_V^i)=0$ for all $i=1,\ldots, \dim X-1$. 
\end{defn}

\begin{exam}
\label{exam:diag}
Let $E$ an elliptic curve over $k$ with an automorphism 
$\rho$ of order 3 and $A:=E^3$.  
The quotient $A/\rho$ (diagonal action)
admits an desingularization $V$ with $K_V=0$. 

There are many embeddings $\iota\,:\, E\hookrightarrow A$ and, 
in particular, every torsion point in $A$
lies on some $\iota(E)$.
 
If $k$ is finite then every point in $V(\bar{k})$
lies on some $\bar{k}$-rational curve in $V$.
Moreover, note that the $E^2/\rho$ (diagonal action) 
is a rational surface.
Hence every point in $V(\bar{k})$ lies in fact on a rational surface 
(defined over $\bar{k}$). 
\end{exam}

\begin{exam}
\label{exam:klein}
Let $C$ be the Klein quartic curve
and  $J:=J_C$ its 
Jacobian. Then the quotient of $J/\sigma$, where $\sigma$ is an
automorphism of order 7, admits a desingularization $V$ which 
is a Calabi-Yau threefold. 
Again, over finite fields, one can show that every algebraic point 
of $V$ lies on a rational curve. 
\end{exam}

\begin{exam}
\label{exam:voisin}
The following varieties have been considered in \cite{voisin}:
Let $S$ be a K3 surface with an involution $\sigma$ and $E$ an elliptic curve
with the standard involution $\tau$. Then a nonsingular model $V$ of 
$E\times S/(\tau\times \sigma)$, is a Calabi-Yau threefold. 
If we choose $S$ and $E$, defined over a finite field, 
so that every algebraic point of $S$ lies on a rational curve, 
then the same property holds for $V$.  
\end{exam}

\begin{conj}
\label{conj:main}
Let $X$ be any smooth projective variety over a finite field $k$.
Assume that $X$ has trivial canonical class and that  $X_{\bar{k}}$ has
trivial algebraic fundamental group. 
Then every algebraic point of $X$ lies on a rational curve $C\subset X$, 
defined over $\bar{k}$.  
\end{conj}

\begin{rem}
Note that if $A$ is a general abelian variety of dimension 
$n\ge 3$ (over $\C$ or over an algebraic closure of $\bar{\mathbb F}_q(x)$) 
and $\sigma$ is the standard involution, then $A/\sigma$
does not contain rational curves, has trivial fundamental group and 
has Kodaira dimension zero (see Remark~\ref{rem:oort}).  
However, the canonical class of a desingularization is nontrivial, for $n\ge 3$.
This also shows that the presence of rational curves is highly 
unstable under deformations. 
\end{rem}

An interesting test of Conjecture~\ref{conj:main} would be the case of a smooth
quintic in $\P^4$. 
For some intriguing connections 
between the counting of points over towers
of finite fields on the Calabi-Yau quintic (the Hasse-Weil zeta function) and 
mirror symmetry we refer to \cite{candelas}.

\bibliographystyle{smfplain}
\bibliography{k3curve}

\end{document}